\documentclass{amsart}
\usepackage{lastpage}
\usepackage{yhmath,verbatim,mathabx}
\usepackage[all]{xy}

\DeclareFontFamily{U}{mathx}{\hyphenchar\font45}
\DeclareFontShape{U}{mathx}{m}{n}{
	<5> <6> <7> <8> <9> <10>
	<10.95> <12> <14.4> <17.28> <20.74> <24.88>
	mathx10
}{}
\DeclareSymbolFont{mathx}{U}{mathx}{m}{n}
\DeclareFontSubstitution{U}{mathx}{m}{n}
\DeclareMathAccent{\widecheck}{0}{mathx}{"71}
\DeclareMathAccent{\wideparen}{0}{mathx}{"75}

\newcommand{\bdis}{\begin{displaymath}}
\newcommand{\edis}{\end{displaymath}}
\newcommand{\be}{\begin{equation}}
\newcommand{\ee}{\end{equation}}
\newcommand{\mbb}{\mathbb}
\newcommand{\mcal}{\mathcal}

\newcommand{\vp}{\varphi}

\newcommand{\zf}{\zeta\left(\frac{1}{2}+it\right)}

\theoremstyle{definition}

\theoremstyle{remark}
\newtheorem{remark}[]{Remark}

\newtheorem*{mydef11}{{\bf Theorem 1}}

\newtheorem*{mydef12}{{\bf Theorem 2}}

\newtheorem*{mydef41}{{\bf Corollary 1}}

\newtheorem*{mydef42}{{\bf Corollary 2}}

\newtheorem*{mydef51}{{\bf Lemma 1}}

\newtheorem*{mydef52}{{\bf Lemma 2}}

\newtheorem*{mydef53}{{\bf Lemma 3}}

\numberwithin{equation}{section}

\begin{document}

\title{Jacob's ladders, crossbreeding and new synergetic formulas for the class of more complicated external parts of $\zeta$-factorization formulas} 

\author{Jan Moser}

\address{Department of Mathematical Analysis and Numerical Mathematics, Comenius University, Mlynska Dolina M105, 842 48 Bratislava, SLOVAKIA}

\email{jan.mozer@fmph.uniba.sk}

\keywords{Riemann zeta-function}

\begin{abstract}
In this paper we obtain new canonical synergetic formula, namely an $\zeta$-analogue of next elementary trigonometric formula. This one describes cooperative interactions between corresponding class of elementary functions and the Riemann's zeta-function on a class of disconnected sets on the critical line. 
\end{abstract}
\maketitle

\section{Introduction}

\subsection{} 

In this paper we obtain new results of the following type: the set of elementary functions 
\be \label{1.1} 
\begin{split} 
& \{t\sin^2t, t\cos^2t,t\cos(2t)\}, \\ 
& t\in [\pi L,\pi L+U],\ U\in(0,\pi/4),\ L\in\mbb{N} 
\end{split} 
\ee  
generates the following synergetic (cooperative) formula 
\be \label{1.2} 
\begin{split}
& \{\alpha_0^{2,1}\tilde{Z}^2(\alpha_1^{2,1})\}\cos^2(\alpha_0^{2,1})-
\{\alpha_0^{1,1}\tilde{Z}^2(\alpha_1^{1,1})\}\sin^2(\alpha_0^{1,1})= \\ 
& = \{\alpha_0^{3,1}\tilde{Z}^2(\alpha_1^{3,1})\}\cos(2\alpha_0^{3,1}),\ L\geq L_0>0, 
\end{split} 
\ee 
with $L_0\in\mbb{N}$ being sufficiently big, where 
\be \label{1.3} 
\tilde{Z}^2(t)=\frac{|\zf|^2}{\omega(t)},\ \omega(t)=\left\{ 1+\mcal{O}\left(\frac{\ln\ln t}{\ln t}\right)\right\}\ln t, 
\ee 
(see \cite{2}, (6.1), (6.7), (7.7), (7.8), (9.1)), next 
\be \label{1.4} 
\begin{split}
& \alpha_0^{1,1},\alpha_0^{2,1},\alpha_0^{3,1}\in (\pi L,\pi L+U), \\ 
& \alpha_1^{1,1},\alpha_1^{2,1},\alpha_1^{3,1}\in (\overset{1}{\wideparen{\pi L}},\overset{1}{\wideparen{\pi L+U}}), 
\end{split}
\ee  
and the segment 
\bdis 
[\overset{1}{\wideparen{\pi L}},\overset{1}{\wideparen{\pi L+U}}]
\edis  
is the first reverse iteration (by means of Jacob's ladder $\vp_1(t)$, see \cite{3}) of the basic segment 
\bdis 
[\pi L,\pi L+U]=[\overset{0}{\wideparen{\pi L}},\overset{0}{\wideparen{\pi L+U}}]. 
\edis 

\subsection{} 

Let us notice that in our theory the following is true: the components of the main $\zeta$-disconnected  set (that is (\ref{1.2}) in our case) 
\be \label{1.5} 
\Delta(\pi L,U,1)=[\pi L,\pi L+U]\bigcup [\overset{1}{\wideparen{\pi L}},\overset{1}{\wideparen{\pi L+U}}] 
\ee  
are separated each from other by the gigantic distance $\rho$: 
\be \label{1.6} 
\begin{split}
& \rho\{[\pi L,\pi L+U]; [\overset{1}{\wideparen{\pi L}},\overset{1}{\wideparen{\pi L+U}}] \}\sim (1-c)\pi(\pi L)\sim \\ 
& \sim\pi (1-c)\frac{L}{\ln L}\to\infty,\ L\to\infty, 
\end{split}
\ee  
($c$ stands for Euler's constant and $\pi(x)$ for the prime-counting function). 

\subsection{} 

Since (see (\ref{1.4})) 
\be \label{1.7} 
\begin{split}
& \frac{\pi L}{\pi L+U}<\frac{\alpha_0^{1,1}}{\alpha_0^{2,1}},\frac{\alpha_0^{3,1}}{\alpha_0^{2,1}}<\frac{\pi L+U}{\pi L} \ \Rightarrow \\ 
& \frac{\alpha_0^{1,1}}{\alpha_0^{2,1}},\frac{\alpha_0^{3,1}}{\alpha_0^{2,1}}\xrightarrow{L\to\infty}1
\end{split}
\ee  
then we have (see (\ref{1.2}), (\ref{1.3}), (\ref{1.7})) the following canonical synergetic formula 
\be \label{1.8} 
\begin{split}
& \left|\zeta\left(\frac 12+i\alpha_1^{2,1}\right)\right|^2\cos^2\alpha_0^{2,1}-
\left|\zeta\left(\frac 12+i\alpha_1^{1,1}\right)\right|^2\sin^2\alpha_0^{1,1}\sim \\ 
& \sim \left|\zeta\left(\frac 12+i\alpha_1^{3,1}\right)\right|^2\cos 2\alpha_0^{3,1},\ 
L\to \infty . 
\end{split}
\ee 

\begin{remark}
Our formula (\ref{1.8}) is: 
\begin{itemize}
	\item[(a)] simple case of formula that is generated by immediate metamorphosis of the main formula, 
	\item[(b)] $\zeta$-analogue of the elementary trigonometric formula 
	\bdis 
	\cos^2x-\sin^2x=\cos 2x , 
	\edis  
	\item[(c)] synergetic one (as well as (\ref{1.2})) since it is generated by interactions between the continuum sets (see our interpretation in \cite{8}) 
	\bdis 
	\left\{\left|\zf\right|^2\right\}, \{t\sin^2t\},\ \{t\cos^2t\},\ \{t\cos 2t\},\ t\geq L_0 
	\edis 
	(some analogue of the classical Belousov-Zhabotiski chemical oscillations), 
	\item[(d)] new type of result simultaneously in the theory of Riemann's zeta-function and in the theory of real continuous functions. 
\end{itemize}
\end{remark} 

\subsection{} 

Finally, we notice the following. 

\begin{remark}
 The formulations of all results and proofs in this paper are based on new notions and methods in the theory of Riemann's zeta-function we have introduced in our series of 47 papers concerning Jacob's ladders. These can be found in arXiv [math.CA] starting with the paper \cite{1}. 
 
 Here we use especially the following notions: Jacob's ladder, $\zeta$-disconnected set that generates the Jacob's ladder, (see \cite{3}), algorithm for generating the $\zeta$-factorization formulas (see \cite{4}), crossbreeding, secondary crossbreeding, exact and asymptotic complete hybrid formula (see \cite{6} -- \cite{9}). Short survey of these notions are listed in papers \cite{5}, \cite{8}. 
\end{remark}

\section{Lemmas} 

By making use of our algorithm for generating $\zeta$-factorization formulas (see \cite{5}, (3.1) -- (3.11), comp. \cite{4}) we obtain the following set of results. 

\subsection{} 

Since 
\bdis 
\begin{split}
& \frac 1U\int_{\pi L}^{\pi L+U}t\sin^2t{\rm d}t= \\ 
& = \frac 14(2\pi L+U)-\frac 12(\pi L+U)\frac{\sin 2U}{2U}+\frac 14\frac{\sin^2U}{U}, 
\end{split}
\edis  
then we obtain the following statement. 

\begin{mydef51}
For the function 
\be \label{2.1} 
\begin{split}
& f_1(t)=t\sin^2t \in \tilde{C}_0[\pi L,\pi L+U],\\ 
& U\in (0,\pi/4),\ L\in\mbb{N}
\end{split}
\ee  
there are vector-valued functions 
\be \label{2.2} 
\begin{split}
& (\alpha_0^{1,k_1},\alpha_1^{1,k_1},\dots,\alpha_{k_1}^{1,k_1},\beta_1^{k_1},\dots,\beta_{k_1}^{k_1}), \\ 
& 1\leq k_1\leq k_0,\ k_1,k_0\in\mbb{N} 
\end{split}
\ee 
(here, we fix arbitrary integer $k_0$) such that the following exact $\zeta$-factorization formula 
\be \label{2.3} 
\begin{split}
& \prod_{r=1}^{k_1}\frac{\tilde{Z}^2(\alpha_r^{1,k_1})}{\tilde{Z}(\beta_r^{k_1})}= 
\frac{1}{\alpha_0^{1,k_1}\sin^2\alpha_0^{1,k_1}}\times \\ 
& \times 
\left\{
\frac 14(2\pi L+U)-\frac 12(\pi L+U)\frac{\sin 2U}{2U}+\frac 14\frac{\sin^2U}{U}
\right\},\ \forall\- L\geq L_0>0 
\end{split}
\ee  
($L_0$ is a sufficiently big one) holds true, where 
\be \label{2.4} 
\begin{split}
& \alpha_r^{1,k_1}=\alpha_r(U,\pi L,k_1;f_1),\ r=0,1,\dots,k_1, \\ 
& \beta_r^{k_1}=\beta_r(U,\pi L,k_1),\ r=1,\dots,k_1, \\ 
& \alpha_0^{1,k_1}\in (\pi L,\pi L+U), 
\alpha_r^{1,k_1}, \beta_r^{k_1}\in 
(\overset{r}{\wideparen{\pi L}},\overset{r}{\wideparen{\pi L+U}}),\ r=1,\dots,k_1, 
\end{split}
\ee  
and the segment 
\bdis 
[\overset{r}{\wideparen{\pi L}},\overset{r}{\wideparen{\pi L+U}}]
\edis  
is the $r$-th reverse iteration by means of the Jacob's ladder, see \cite{3}, of the basic segment 
\bdis 
[\pi L,\pi L+U]=[\overset{0}{\wideparen{\pi L}},\overset{0}{\wideparen{\pi L+U}}] . 
\edis 
\end{mydef51} 

\subsection{} 

Since 
\bdis 
\begin{split}
	& \frac 1U\int_{\pi L}^{\pi L+U}t\cos^2t{\rm d}t= \\ 
	& = \frac 14(2\pi L+U)+\frac 12(\pi L+U)\frac{\sin 2U}{2U}-\frac 14\frac{\sin^2U}{U}, 
\end{split}
\edis  
then we obtain the following statement. 

\begin{mydef52}
	For the function 
	\be \label{2.5} 
	\begin{split}
		& f_2(t)=t\cos^2t \in \tilde{C}_0[\pi L,\pi L+U],\\ 
		& U\in (0,\pi/4) 
	\end{split}
	\ee  
	there are vector-valued functions 
	\be \label{2.6} 
	\begin{split}
		& (\alpha_0^{2,k_2},\alpha_1^{2,k_2},\dots,\alpha_{k_2}^{2,k_2},\beta_1^{k_2},\dots,\beta_{k_2}^{k_2}), \\ 
		& 1\leq k_2\leq k_0,\ k_2\in\mbb{N} 
	\end{split}
	\ee 
	such that the following exact $\zeta$-factorization formula 
	\be \label{2.7} 
	\begin{split}
		& \prod_{r=1}^{k_2}\frac{\tilde{Z}^2(\alpha_r^{2,k_2})}{\tilde{Z}(\beta_r^{k_2})}= 
		\frac{1}{\alpha_0^{2,k_2}\cos^2\alpha_0^{2,k_2}}\times \\ 
		& \times 
		\left\{
		\frac 14(2\pi L+U)+\frac 12(\pi L+U)\frac{\sin 2U}{2U}-\frac 14\frac{\sin^2U}{U}
		\right\},\ \forall\- L\geq L_0>0 
	\end{split}
	\ee  
	holds true, where 
	\be \label{2.8} 
	\begin{split}
		& \alpha_r^{2,k_2}=\alpha_r(U,\pi L,k_2;f_2),\ r=0,1,\dots,k_2, \\ 
		& \beta_r^{k_2}=\beta_r(U,\pi L,k_2),\ r=1,\dots,k_2, \\ 
		& \alpha_0^{2,k_2}\in (\pi L,\pi L+U), 
		\alpha_r^{2,k_2}, \beta_r^{k_2}\in 
		(\overset{r}{\wideparen{\pi L}},\overset{r}{\wideparen{\pi L+U}}),\ r=1,\dots,k_2 . 
	\end{split}
	\ee  
\end{mydef52}  

\subsection{} 

Since 
\bdis 
\begin{split}
	& \frac 1U\int_{\pi L}^{\pi L+U}t\cos 2t{\rm d}t= 
(\pi L+U)\frac{\sin 2U}{2U}-\frac 12\frac{\sin^2U}{U}, 
\end{split}
\edis  
then we obtain the following statement. 

\begin{mydef53}
	For the function 
	\be \label{2.9} 
	\begin{split}
		& f_3(t)=t\cos 2t \in \tilde{C}_0[\pi L,\pi L+U],\\ 
		& U\in (0,\pi/4) 
	\end{split}
	\ee  
	there are vector-valued functions 
	\be \label{2.10} 
	\begin{split}
		& (\alpha_0^{3,k_3},\alpha_1^{3,k_3},\dots,\alpha_{k_3}^{3,k_3},\beta_1^{k_3},\dots,\beta_{k_3}^{k_3}), \\ 
		& 1\leq k_3\leq k_0,\ k_3\in\mbb{N} 
	\end{split}
	\ee 
	such that the following exact $\zeta$-factorization formula 
	\be \label{2.11} 
	\begin{split}
		& \prod_{r=1}^{k_3}\frac{\tilde{Z}^2(\alpha_r^{3,k_3})}{\tilde{Z}(\beta_r^{k_3})}= \\ 
		& = 
		\frac{1}{\alpha_0^{3,k_3}\cos (2\alpha_0^{3,k_3})}\times  
		\left\{
		(\pi L+U)\frac{\sin 2U}{2U}-\frac 12\frac{\sin^2U}{U}
		\right\},\ \forall\- L\geq L_0>0 
	\end{split}
	\ee  
	holds true, where 
	\be \label{2.12} 
	\begin{split}
		& \alpha_r^{3,k_3}=\alpha_r(U,\pi L,k_3;f_3),\ r=0,1,\dots,k_3, \\ 
		& \beta_r^{k_3}=\beta_r(U,\pi L,k_3),\ r=1,\dots,k_3, \\ 
		& \alpha_0^{3,k_3}\in (\pi L,\pi L+U), 
		\alpha_r^{3,k_3}, \beta_r^{k_3}\in 
		(\overset{r}{\wideparen{\pi L}},\overset{r}{\wideparen{\pi L+U}}),\ r=1,\dots,k_3 . 
	\end{split}
	\ee  
\end{mydef53}  

\section{Exact complete hybrid formula} 

\subsection{} 

We start with the following. 

\begin{remark}
Our description of the operation of crossbreeding (see \cite{6} and \cite{8}, subsection 3.3) contains the following expression: 

\dots that is: after finite number of eliminations of the external functions 
\bdis 
E_m(U,T),\ m=1,\dots,M,\dots 
\edis  
However, this is not exact. The exact phrase is as follows: 

\dots that is: after finite number of eliminations of the variables $U,T$ from the set of external functions \dots 

We shall call these variables as external ones. 
\end{remark} 

\subsection{} 

Now, we make the crossbreeding on the set 
\be \label{3.1} 
\{ (2.3), (2.7), (2.11)\},\ U\in (0,\pi/4) 
\ee  
of exact $\zeta$-factorization formulas. 

\begin{remark}
In the case (\ref{3.1}) we see that the corresponding external functions contain the pair $\pi L, U$ of external variables, comp. Remark 3. 
\end{remark}

First, elimination of the block $\{\dots\}$ (see (\ref{2.11})) from (\ref{2.3}), (\ref{2.11}) gives 
\be \label{3.2} 
\begin{split}
& \alpha_0^{1,k_1}\sin^2\alpha_0^{1,k_1}\prod_{r=1}^{k_1}
\frac{\tilde{Z}^2(\alpha_r^{1,k_1})}{\tilde{Z}^2(\beta_r^{k_1})}+\\ 
& + \frac 12\alpha_0^{3,k_3}\cos(2\alpha_0^{3,k_3})\prod_{r=1}^{k_3}
\frac{\tilde{Z}^2(\alpha_r^{3,k_3})}{\tilde{Z}^2(\beta_r^{k_3})}=\frac 12\left(\pi L+\frac U2\right),
\end{split}
\ee 
and socondly, (\ref{2.3}) and (\ref{2.7}) imply 
\be \label{3.3} 
\begin{split}
	& \alpha_0^{2,k_2}\cos^2\alpha_0^{2,k_2}\prod_{r=1}^{k_2}
	\frac{\tilde{Z}^2(\alpha_r^{2,k_2})}{\tilde{Z}^2(\beta_r^{k_2})}+\\ 
	& + \alpha_0^{1,k_1}\sin^2\alpha_0^{1,k_1}\prod_{r=1}^{k_1}
	\frac{\tilde{Z}^2(\alpha_r^{1,k_1})}{\tilde{Z}^2(\beta_r^{k_1})}=\pi L+\frac U2. 
\end{split}
\ee  
Finally, we obtain from (\ref{3.2}) and (\ref{3.3}) the following 

\begin{mydef11}
The set 
\be \label{3.4} 
\begin{split} 
	& \{t\sin^2t, t\cos^2t,t\cos(2t)\}, \\ 
	& t\in [\pi L,\pi L+U],\ U\in(0,\pi/4),\ L\in\mbb{N} 
\end{split} 
\ee 
of elementary functions generates the following exact complete hybrid formula 
\be \label{3.5} 
\begin{split}
& \left\{
\alpha_0^{2,k_2}\prod_{r=1}^{k_2}
\frac{\tilde{Z}^2(\alpha_r^{2,k_2})}{\tilde{Z}^2(\beta_r^{k_2})}
\right\}\cos^2\alpha_0^{2,k_2}- 
\left\{
\alpha_0^{1,k_1}\prod_{r=1}^{k_1}
\frac{\tilde{Z}^2(\alpha_r^{1,k_1})}{\tilde{Z}^2(\beta_r^{k_1})}
\right\}\sin^2\alpha_0^{1,k_1}= \\ 
& = 
\left\{
\alpha_0^{3,k_3}\prod_{r=1}^{k_3}
\frac{\tilde{Z}^2(\alpha_r^{3,k_3})}{\tilde{Z}^2(\beta_r^{k_3})}
\right\}\cos(2\alpha_0^{3,k_3}), \\ 
& \forall\- L\geq L_0>0,\ 1\leq k_1,k_2,k_3\leq k_0, 
\end{split}
\ee  
(wee fix arbitrary $k_0\in\mbb{N}$ and $L_0$ is a sufficiently big one), where 
\be \label{3.6} 
\begin{split}
& \alpha_0^{1,k_1},\alpha_0^{2,k_2}, \alpha_0^{3,k_3}\in (\pi L,\pi L+U), \\ 
& \alpha_r^{1,k_1},\beta_r^{k_1}\in (\overset{r}{\wideparen{\pi L}},\overset{r}{\wideparen{\pi L+U}}),\ r=1,\dots,k_1, \\ 
& \alpha_r^{2,k_2},\beta_r^{k_2}\in (\overset{r}{\wideparen{\pi L}},\overset{r}{\wideparen{\pi L+U}}),\ r=1,\dots,k_2, \\ 
& \alpha_r^{3,k_3},\beta_r^{k_3}\in (\overset{r}{\wideparen{\pi L}},\overset{r}{\wideparen{\pi L+U}}),\ r=1,\dots,k_3,
\end{split}
\ee 
i.e. 
\bdis 
\begin{split}
& r=0,1,\dots,\bar{k}:\ \alpha_r^{1,k_1},\alpha_r^{2,k_2},\alpha_r^{3,k_3}, \beta_r^{k_1}, \beta_r^{k_2}, \beta_r^{k_3}\in  \\ 
& \in \Delta(U,\pi L,\bar{k})=\bigcup_{r=0}^{\bar{k}}
[\overset{r}{\wideparen{\pi L}},\overset{r}{\wideparen{\pi L+U}}],\ \bar{k}=\max\{ k_1,k_2,k_3\},\ 
\bar{k}\leq k_0, 
\end{split}
\edis  
and 
\be \label{3.7} 
\Delta(U,\pi L,\bar{k})\subset \Delta(U,\pi L,k_0)=
\bigcup_{r=0}^{k_0}
[\overset{r}{\wideparen{\pi L}},\overset{r}{\wideparen{\pi L+U}}], 
\ee  
where the last $\zeta$-disconnected set is basic one (for every fixed $k_0\in \mbb{N}$). 
\end{mydef11} 

\begin{remark}
It is true in our theory (see \cite{3}): consecutive components of the basic disconnected set are separated each from other by gigantic distances $\rho$: 
\be \label{3.8} 
\begin{split}
& \rho
\left\{
[\overset{r}{\wideparen{\pi L}},\overset{r}{\wideparen{\pi L+U}}], 
[\overset{r+1}{\wideparen{\pi L}},\overset{r+1}{\wideparen{\pi L+U}}]
\right\}\sim (1-c)\pi(\pi L)\sim \\ 
& \sim 
\pi \cdot (1-c)\frac{L}{\ln L}\xrightarrow{L\to\infty}\infty,\ r=0,1,\dots,k_0-1,  
\end{split}
\ee  
($c$ is the Euler's constant and $\pi(x)$ stands for the prime-counting function). 
\end{remark} 

\begin{remark}
By our interpretation given in the paper \cite{8} the formula (\ref{3.5}) is the synergetic (cooperative) one in the following sense: it is the result of interactions between the following continuum sets 
\be \label{3.9} 
\left\{
\left|\zf\right|^2
\right\}, \{ t\sin^2t\},\ \{ t\cos^2t\}, \{ t\cos 2t\},\ t\geq L_0 
\ee  
and these interactions are excited by the Jacob's ladder $\vp_1(t)$. We call these interactions (see \cite{9}) as the $\zeta$-chemical reaction between sets (\ref{3.9}). 
\end{remark} 

\begin{remark}
The result of above mentioned $\zeta$-chemical reactions (the $\zeta$-chemical compound) is our synergetic formula (\ref{3.5}). This interpretation represents a $\zeta$-analogue of the classical Belousov-Zhabotiski chemical oscillations (see our paper \cite{8} as the starting point in this direction). 
\end{remark} 

\section{Immediate metamorphosis of the formula (\ref{3.5}) into asymptotic secondary complete hybrid formula} 

\subsection{} 

If we rewrite the formula (\ref{3.5}) in the form 
\bdis 
\{\dots\}\cos^2\alpha_0^{2,k_2}=\{\dots\}\sin^2\alpha_0^{1,k_1}+\{\dots\}\cos(2\alpha_0^{3,k_3})
\edis 
and use (\ref{1.3}), (\ref{1.7}) and some small algebra (comp. \cite{8}, Section 8.2), then we obtain the following 

\begin{mydef41}
\be \label{4.1} 
\begin{split}
& \left\{
\prod_{r=1}^{k_2}
\frac{\left|\zeta\left(\frac 12+i\alpha_r^{2,k_2}\right)\right|^2}
{\left|\zeta\left(\frac 12+i\beta_r^{k_2}\right)\right|^2}
\right\}\cos^2\alpha_0^{2,k_2}-
\left\{
\prod_{r=1}^{k_1}
\frac{\left|\zeta\left(\frac 12+i\alpha_r^{1,k_1}\right)\right|^2}
{\left|\zeta\left(\frac 12+i\beta_r^{k_1}\right)\right|^2}
\right\}\sin^2\alpha_0^{1,k_1}\sim \\ 
& \sim 
\left\{
\prod_{r=1}^{k_3}
\frac{\left|\zeta\left(\frac 12+i\alpha_r^{3,k_3}\right)\right|^2}
{\left|\zeta\left(\frac 12+i\beta_r^{k_3}\right)\right|^2}
\right\}\cos (2\alpha_0^{3,k_3}),\ L\to\infty. 
\end{split}
\ee 
\end{mydef41} 

\subsection{} 

Now, in the case 
\be \label{4.2} 
k_1=k_2=k_3=k;\ 1\leq k\leq k_0 
\ee  
we obtain the following 

\begin{mydef42}
\be \label{4.3} 
\begin{split}
& \left\{
\prod_{r=1}^k\left|\zeta\left(\frac 12+i\alpha_r^{2,k}\right)\right|^2
\right\}\cos^2\alpha_0^{2,k}-
\left\{
\prod_{r=1}^k\left|\zeta\left(\frac 12+i\alpha_r^{1,k}\right)\right|^2
\right\}\sin^2\alpha_0^{1,k}\sim \\ 
& \sim 
\left\{
\prod_{r=1}^k\left|\zeta\left(\frac 12+i\alpha_r^{3,k}\right)\right|^2
\right\}\cos (2\alpha_0^{3,k}),\ L\to\infty. 
\end{split}
\ee  
\end{mydef42} 

\begin{remark}
Formula (\ref{4.3}) expresses the result of immediate metamorphosis of the asymptotic complete hybrid formula (\ref{4.1}) into asymptotic secondary complete hybrid formula in the case (\ref{4.2}). 
\end{remark} 

\begin{remark}
The case $k=1$ in (\ref{4.3}) gives the formula (\ref{1.8}) that has been used in Introduction to inform about the content of this paper. 
\end{remark} 

\section{Secondary exact complete hybrid formula} 

We choose the following exact complete hybrid formula (see \cite{9}, (3.7)) 
\be \label{5.1} 
\begin{split}
	& (1+\Delta_4)^{1/\Delta_4}(\alpha_0^{4,k_4}-\pi L)
	\left\{
	\prod_{r=1}^{k_4}\frac{\tilde{Z}^2(\alpha_r^{4,k_4})}{\tilde{Z}^2(\beta_r^{k_4})}
	\right\}^{1/\Delta_4}= \\ 
	& = (1+\Delta_5)^{1/\Delta_5}(\alpha_0^{5,k_5}-\pi L)
	\left\{
	\prod_{r=1}^{k_5}\frac{\tilde{Z}^2(\alpha_r^{5,k_5})}{\tilde{Z}^2(\beta_r^{k_5})}
	\right\}^{1/\Delta_5}, \\ 
	& \forall\- L\geq L_0>0,\ \Delta_4,\Delta_5>0,\ \Delta_5\not=\Delta_4,\ 1\leq k_4,k_5\leq k_0. 
\end{split}
\ee  
Now, we make the use of operation of secondary crossbreeding (see \cite{8}) on the set 
\bdis 
\{(3.5), (5.1)\} 
\edis 
as follows. First of all, we put 
\bdis 
k_4=k_5=k;\ 1\leq k\leq k_0 
\edis  
in the formula (\ref{5.1}) that gives the result 
\be \label{5.2} 
\begin{split}
	& \prod_{r=1}^k\tilde{Z}^2(\beta_r^k)=\\ 
	& = 
	\left[\frac{(1+\Delta_4)^{1/\Delta_4}}{(1+\Delta_5)^{1/\Delta_5}}\right]^{\frac{\Delta_5\Delta_4}{\Delta_5-\Delta_4}}
	\left( \frac{\alpha_0^{4,k}-\pi L}{\alpha_0^{5,k}-\pi L}\right)^{\frac{\Delta_5\Delta_4}{\Delta_5-\Delta_4}}
	\left\{
	\prod_{r=1}^k\tilde{Z}^2(\alpha_r^{4,k})
	\right\}^{\frac{\Delta_5}{\Delta_5-\Delta_4}}\times \\ 
	& \times 
	\left\{
	\prod_{r=1}^k\tilde{Z}^2(\alpha_r^{5,k})
	\right\}^{-\frac{\Delta_4}{\Delta_5-\Delta_4}}. 
\end{split}
\ee  
For the second, we put consecutively 
\bdis 
k=k_1,k_2,k_3 
\edis  
in (\ref{5.2}) and the corresponding results we substitute into the formula (\ref{3.5}). The final result is expressed by the following 

\begin{mydef12}
The two sets of elementary functions 
\bdis 
\begin{split}
& \{t\sin^2t,t\cos^2t,t\cos 2t\},\ \{ (t-\pi L)^{\Delta_4},(t-\pi L)^{\Delta_5}\}, \\ 
& t\in [\pi L,\pi L+U],\ U\in (0,\pi/4), \Delta_4,\Delta_5>0, \Delta_4\not=\Delta_5 
\end{split}
\edis  
generate  the following secondary exact complete hybrid formula 
\be \label{5.3} 
\begin{split}
& \alpha_0^{2,k_2}
\frac
{\prod_{r=1}^{k_2}\tilde{Z}^2(\alpha_r^{2,k_2})[\tilde{Z}^2(\alpha_r^{5,k_2})]^{\frac{\Delta_4}{\Delta_5-\Delta_4}}[\tilde{Z}^2(\alpha_r^{4,k_2})]^{-\frac{\Delta_5}{\Delta_5-\Delta_4}}}
{\left( \frac{\alpha_0^{4,k_2}-\pi L}{\alpha_0^{5,k_2}-\pi L}\right)^{\frac{\Delta_5\Delta_4}{\Delta_5-\Delta_4}}}\cos^2\alpha_0^{2,k_2}
- \\ 
& - 
\alpha_0^{1,k_1}
\frac
{\prod_{r=1}^{k_1}\tilde{Z}^2(\alpha_r^{1,k_1})[\tilde{Z}^2(\alpha_r^{5,k_1})]^{\frac{\Delta_4}{\Delta_5-\Delta_4}}[\tilde{Z}^2(\alpha_r^{4,k_1})]^{-\frac{\Delta_5}{\Delta_5-\Delta_4}}}
{\left( \frac{\alpha_0^{4,k_1}-\pi L}{\alpha_0^{5,k_1}-\pi L}\right)^{\frac{\Delta_5\Delta_4}{\Delta_5-\Delta_4}}}\sin^2\alpha_0^{1,k_1}= \\ 
& = 
\alpha_0^{3,k_3}
\frac
{\prod_{r=1}^{k_3}\tilde{Z}^2(\alpha_r^{3,k_3})[\tilde{Z}^2(\alpha_r^{5,k_3})]^{\frac{\Delta_4}{\Delta_5-\Delta_4}}[\tilde{Z}^2(\alpha_r^{4,k_3})]^{-\frac{\Delta_5}{\Delta_5-\Delta_4}}}
{\left( \frac{\alpha_0^{4,k_3}-\pi L}{\alpha_0^{5,k_3}-\pi L}\right)^{\frac{\Delta_5\Delta_4}{\Delta_5-\Delta_4}}}\cos (2\alpha_0^{3,k_3}), \\ 
& \forall\- L\geq L_0, 1\leq k_1,k_2,k_3\leq k_0. 
\end{split}
\ee 
\end{mydef12} 

\begin{remark}
Let us notice explicitly that  two complicated types of $\zeta$-modulation (of amplitude and also phase) of the elementary trigonometric formula 
\bdis 
\cos^2x-\sin^2x=\cos 2x 
\edis 
are expressed by the synergetic formulae (\ref{3.5}), (\ref{5.3}). 
\end{remark}

I would like to thank Michal Demetrian for his moral support of my study of Jacob's ladders.

\end{document}